\documentclass[11pt, reqno]{amsart}
\usepackage{natbib}
\usepackage{latexcad}

\def\adjustpagesize#1#2#3#4{
    \addtolength{\textwidth}{#1mm}\addtolength{\hoffset}{#2mm}
    \addtolength{\textheight}{#3mm}\addtolength{\voffset}{#4mm}
}

\theoremstyle{definition}
\newtheorem{defi}{Definition}[section]
\newtheorem{exam}[defi]{Example}
\newtheorem{assu}[defi]{Assumption}

\theoremstyle{remark}

\theoremstyle{plain}
\newtheorem{theo}[defi]{Theorem}
\newtheorem{lemm}[defi]{Lemma}

\newtheorem{prop}[defi]{Proposition}

\newcommand{\theolbl}[1]{\label{theo:#1}}
\newcommand{\lemmlbl}[1]{\label{lemm:#1}}

\newcommand{\proplbl}[1]{\label{prop:#1}}

\newcommand{\assulbl}[1]{\label{assu:#1}}

\newcommand{\theoref}[1]{\ref{theo:#1}}
\newcommand{\lemmref}[1]{\ref{lemm:#1}}

\newcommand{\theofref}[1]{Theorem \ref{theo:#1}}
\newcommand{\lemmfref}[1]{Lemma \ref{lemm:#1}}

\newcommand{\propfref}[1]{Proposition \ref{prop:#1}}

\newcommand{\assufref}[1]{Assumption \ref{assu:#1}}

\def\var#1{\mathrm{var}\,#1}

\adjustpagesize{20}{-10}{30}{-15} 

\title[Loops of Bol-Moufang type]{Loops of Bol-Moufang type with a subgroup of index two}

\author[M.~K.~Kinyon]{Michael~K.~Kinyon}
\address{Department of Mathematical Sciences\\
Indiana University South Bend \\
South Bend, Indiana 46634\\
U.S.A}
\email{mkinyon@iusb.edu}
\urladdr{http://mypage.iusb.edu/\symbol{126}mkinyon}

\author[J.~D.~Phillips]{J.~D.~Phillips}
\address{Department of Mathematics \& Computer Science\\
Wabash College\\
Crawfordsville, Indiana 47933\\
U.S.A.}
\email{phillipj@wabash.edu}

\author[P.~Vojt\v{e}chovsk\'y]{Petr Vojt\v{e}chovsk\'y}
\address{Department of Mathematics\\
University of Denver\\
2360 S Gaylord St\\
Denver, Colorado 80208\\
U.S.A.}
\email{petr@math.du.edu}
\urladdr{http://www.math.du.edu/\symbol{126}petr}

\begin{document}

\begin{abstract} We describe all constructions for loops of Bol-Moufang type
analogous to the Chein construction $M(G,*,g_0)$ for Moufang loops.
\end{abstract}



\maketitle

\section{Introduction}

Due to the specialized nature of this paper we assume that the reader is
already familiar with the theory of quasigroups and loops. We therefore omit
basic definitions and results (see \cite{Bruck}, \cite{Pflugfelder}).

In a sense, a nonassociative loop is closest to a group when it contains a
subgroup of index two. Such loops proved useful in the study of Moufang
loops, and it is our opinion that they will also prove useful in the study of
other varieties of loops.

Here is the well-known construction of Moufang loops with a subgroup of index
two:

\begin{theo}[Chein \cite{Chein}]\theolbl{Chein}
Let $G$ be a group, $g_0\in Z(G)$, and $*$ an involutory
antiautomorphism of $G$ such that $g_0^*=g_0$, $gg^*\in Z(G)$ for every $g\in
G$. For an indeterminate $u$, define multiplication $\circ$ on $G\cup Gu$ by
\begin{equation}\label{Eq:Chein}
    g\circ h = gh,\quad g\circ (hu) = (hg)u,\quad
    gu\circ h = (gh^*)u,\quad gu\circ hu = g_0h^*g,
\end{equation}
where $g$, $h\in G$. Then $L=(G\cup Gu,\circ)$ is a Moufang loop. Moreover,
$L$ is associative if and only if $G$ is commutative.
\end{theo}

It has been shown in \cite{VojtechovskyCMUC1} that \eqref{Eq:Chein} is the
only construction of its kind for Moufang loops. (This statement will be
clarified later.) In \cite{VojtechovskyCMUC2}, all constructions similar to
\eqref{Eq:Chein} were determined for Bol loops.

%FIGURE
\setlength{\unitlength}{1.00mm}
\begin{figure}[ht]\begin{center}\input{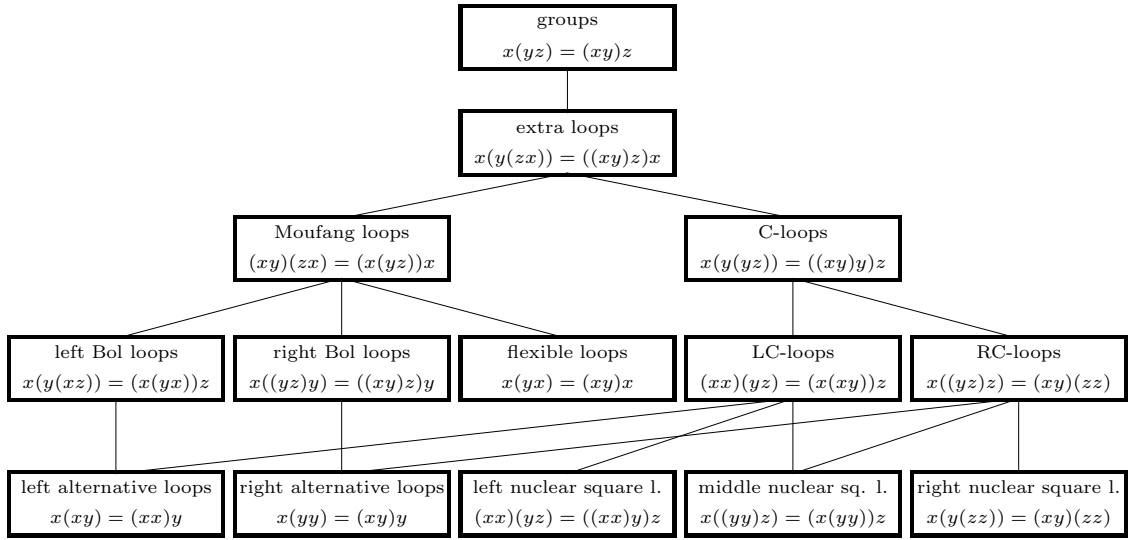}\end{center}
\caption{The varieties of loops of Bol-Moufang type.}\label{Fg:BM}
\end{figure}

The purpose of this paper is to give a complete list of all constructions
similar to \eqref{Eq:Chein} for all loops of Bol-Moufang type. A groupoid
identity is of \emph{Bol-Moufang type} if it has three distinct variables,
two of the variables occur once on each side, the third variable occurs twice
on each side, and the variables occur in the same order on both sides. A loop
is of \emph{Bol-Moufang type} if it belongs to a variety of loops defined by
a single identity of Bol-Moufang type. Figure \ref{Fg:BM} shows all varieties
of loops of Bol-Moufang type and all inclusions among them (cf.
\cite{Fenyves}, \cite{PhillipsVojtechovsky}). Some varieties of Figure
\ref{Fg:BM} can be defined equivalently by other identities of Bol-Moufang
type. For instance, Moufang loops are equivalently defined by the identity
$x(y(xz))=((xy)x)z$. See \cite{PhillipsVojtechovsky} for all such
equivalences. Furthermore, although some defining identities of Figure
\ref{Fg:BM} do not appear to be of Bol-Moufang type, they are in fact
equivalent to some Bol-Moufang identity. For instance, the flexible law
$x(yx)=(xy)x$ is equivalent to the Bol-Moufang identity $(x(yx))z = ((xy)x)z$
in any variety of loops.

As we shall see, the computational complexity of our programme is
overwhelming (for humans). We therefore first carefully define what we mean
by a construction similar to \eqref{Eq:Chein} (see Section \ref{Sc:Mult}),
and then identify situations in which two given constructions are ``the
same'' (see Sections \ref{Sc:Red}, \ref{Sc:Same}, \ref{Sc:Loops}). Upon
showing which constructions yield loops, we work out one construction by hand
(see Section \ref{Sc:dBJ}), and then switch to a computer search, described
in Section \ref{Sc:Search}. The results of the computer search are summarized
in Section \ref{Sc:Results}.

\section{Similar Constructions}\label{Sc:Mult}

Throughout the paper, we assume that $G$ is a finite group, $g_0\in Z(G)$,
and $*$ is an involutory automorphism of $G$ such that $g_0^*=g_0$ and
$gg^*\in Z(G)$ for every $g\in G$.

The following property of $*$ will be used without reference:

\begin{lemm} Let $G$ be a group and $*:G\to G$ an involutory map
such that $gg^*\in Z(G)$ for every $g\in G$. Then $g^*g=gg^*\in Z(G)$ for
every $g\in G$.
\end{lemm}
\begin{proof}
For $g\in G$, we have $g^*g = g^*(g^*)^* \in Z(G)$. Then $(g^*g)g^* =
g^*(g^*g)$, and $gg^*=g^*g$ follows upon cancelling $g^*$ on the left.
\end{proof}

Consider the following eight bijections of $G\times G$:
\begin{displaymath}
    \begin{array}{llll}
        \theta_{xy}(g,h) = (g,h), &\theta_{xy^*}(g,h) = (g,h^*),
        &\theta_{x^*y}(g,h) = (g^*,h), &\theta_{x^*y^*}(g,h) = (g^*,h^*),\\
        \theta_{yx}(g,h) = (h,g), &\theta_{yx^*}(g,h) = (h,g^*),
        &\theta_{y^*x}(g,h) = (h^*,g), &\theta_{y^*x^*}(g,h) = (h^*,g^*).
    \end{array}
\end{displaymath}
They form a group $\Theta$ under composition, isomorphic to the dihedral
group $D_8$. It is generated by $\{\theta_{yx},\theta_{xy^*}\}$, say. Let
$\Theta_0$ be the group generated by $\Theta$ and $\theta_{g_0}$, where
$\theta_{g_0}(g,h)=(g_0g,h)$.

Let $\Delta:G\times G\to G$ be the evaluation map $\Delta(g,h)=gh$, and $u$
an indeterminate. Given $\alpha$, $\beta$, $\gamma$, $\delta\in \Theta_0$,
define multiplication $\circ$ on $G\cup Gu$ by
\begin{displaymath}
    g\circ h = \Delta \alpha(g,h),\quad g\circ hu = (\Delta\beta(g,h))u,\quad
    gu\circ h = (\Delta\gamma(g,h))u,\quad gu\circ hu = \Delta\delta(g,h),
\end{displaymath}
where $g$, $h\in G$. The resulting groupoid $(G\cup Gu,\circ)$ will be
denoted by
\begin{displaymath}
    Q(G,*,g_0,\alpha,\beta,\gamma,\delta),
\end{displaymath}
or by $Q(G,\alpha,\beta,\gamma,\delta)$, when $g_0$, $*$ are known from the
context or if they are not important. It is easy to check that
$Q(G,*,g_0,\alpha,\beta,\gamma,\delta)$ is a quasigroup.

We also define
\begin{displaymath}
    \mathcal Q(G,*,g_0) = \{Q(G,*,g_0,\alpha,\beta,\gamma,\delta);\;
        \alpha,\,\beta,\,\gamma,\,\delta\in \Theta_0\},
\end{displaymath}
and
\begin{displaymath}
    \mathcal Q(G) = \bigcup_{*,g_0} \mathcal Q(G,*,g_0),
\end{displaymath}
where the union is taken over all involutory antiautomorphisms $*$ satisfying
$gg^*\in Z(G)$ for every $g\in G$, and over all elements $g_0$ such that
$g_0^*=g_0\in Z(G)$. By definition, we call elements of $\mathcal Q(G)$
\emph{quasigroups obtained from $G$ by a construction similar to}
\eqref{Eq:Chein}.

\section{Reductions}\label{Sc:Red}

\noindent The goal of this section is to show that one does not have to take
all elements of $\Theta_0$ into consideration in order to determine $\mathcal
Q(G,*,g_0)$.

Note that $g_0^n = (g_0^n)^*\in Z(G)$ for every integer $n$. Therefore
\begin{equation}\label{Eq:g0}
    g_0^n\Delta\theta_0(g,h) = \Delta\theta_{g_0}^n\theta_0(g,h) =
    \Delta\theta_0\theta_{g_0}^n(g,h)
\end{equation}
for every $\theta_0\in\Theta_0$ and every $g$, $h\in G$.

\begin{lemm}
For every integer $n$, the quasigroup
$Q(G,\theta_{g_0}^n\alpha,\theta_{g_0}^n\beta,\theta_{g_0}^n\gamma,
\theta_{g_0}^n\delta)$ is isomorphic to $Q(G,\alpha,\beta,\gamma,\delta)$.
\end{lemm}
\begin{proof} We use \eqref{Eq:g0} freely in this proof. Let $t=g_0^n$.
Denote by $\circ$ the multiplication in $Q(G,\alpha,\beta,\gamma,\delta)$,
and by $\bullet$ the multiplication in $Q(G,\theta_{g_0}^n\alpha$,
$\theta_{g_0}^n\beta$, $\theta_{g_0}^n\gamma$, $\theta_{g_0}^n\delta)$. Let
$f$ be the bijection of $G\cup Gu$ defined by $g\mapsto t^{-1}g$, $gu\mapsto
(t^{-1}g)u$, for $g\in G$. Then for $g$, $h\in G$, we have
\begin{align*}
    &f(g\circ h) = t^{-1}\Delta\alpha(g,h) =
    t\Delta\alpha(t^{-1}g,t^{-1}h) = t^{-1}g\bullet t^{-1}h = f(g)\bullet
    f(h),\\
    &f(g\circ hu) =  t^{-1}\Delta\beta(g,h)u =
    t\Delta\beta(t^{-1}g,t^{-1}h)u = t^{-1}g\bullet (t^{-1}h)u
    = f(g)\bullet f(hu),
\end{align*}
and similarly for $\gamma$, $\delta$. Hence $f$ is the desired isomorphism.
\end{proof}

Therefore, if we only count the quasigroups in $\mathcal Q(G,*,g_0)$ up to
isomorphism, we can assume that $\mathcal Q(G,*,g_0) =
\{Q(G,*,g_0,\alpha,\beta,\gamma,\delta);\; \alpha\in\Theta$, and $\beta$,
$\gamma$, $\delta$ are of the form $\theta\theta_{g_0}^n$ for some
$n\in\mathbb Z$ and $\theta\in\Theta\}$.

Given a groupoid $(A,\cdot)$, the \emph{opposite groupoid}
$(A,\cdot^{\mathrm{op}})$ is defined by $x\cdot^{\mathrm{op}} y = y\cdot x$.

\begin{lemm} The quasigroups $Q(G,\alpha,\beta,\gamma,\delta)$
and $Q(G,\theta_{yx}\alpha$, $\theta_{yx}\gamma$, $\theta_{yx}\beta$,
$\theta_{yx}\delta)$ are opposite to each other.
\end{lemm}
\begin{proof}
Let $\circ$ denote the multiplication in  $Q(G,\alpha,\beta,\gamma,\delta)$,
and $\bullet$ the multiplication in $Q(G,\theta_{yx}\alpha$,
$\theta_{yx}\gamma$, $\theta_{yx}\beta$, $\theta_{yx}\delta\}$. For $g$,
$h\in G$ we have
\begin{align*}
    &g\circ h = \Delta\alpha(g,h) = \Delta\theta_{yx}\alpha(h,g)
        = h\bullet g,\\
    &g\circ hu = \Delta\beta(g,h)u = \Delta\theta_{yx}\beta(h,g) u
        = hu\bullet g,\\
    &gu\circ h = \Delta\gamma(g,h)u = \Delta\theta_{yx}\gamma(h,g) u
        = h\bullet gu,\\
    &gu\circ hu = \Delta\delta(g,h) = \Delta\theta_{yx}\delta(h,g)
        = hu\bullet gu.
\end{align*}
\end{proof}

Therefore, if we only count the quasigroups in $\mathcal Q(G,*,g_0)$ up to
isomorphism and opposites, we can assume that $\mathcal Q(G,*,g_0) =
\{Q(G,*,g_0,\alpha,\beta,\gamma,\delta);\; \alpha\in\{\theta_{xy}$,
$\theta_{xy^*}$, $\theta_{x^*y}$, $\theta_{x^*y^*}\}$, and $\beta$, $\gamma$,
$\delta$ are of the form $\theta\theta_{g_0}^n$ for some $n\in\mathbb Z$ and
$\theta\in\Theta\}$.

\begin{assu}\assulbl{loop} From now on we assume that $\alpha\in\{\theta_{xy}$,
$\theta_{xy^*}$, $\theta_{x^*y}$, $\theta_{x^*y^*}\}$, and that $\beta$,
$\gamma$, $\delta$ are of the form $\theta\theta_{g_0}^n$ for some
$n\in\mathbb Z$ and $\theta\in\Theta$.
\end{assu}

\section{When $*$ is identical on $G$}\label{Sc:Same}

\noindent Assume for a while that $g=g^*$ for every $g\in G$. Then $gh =
(gh)^* = h^*g^* = hg$ shows that $G$ is commutative. In particular, $\Theta =
\{\theta_{xy}\}$, and $\Theta_0 = \bigcup_n \theta_{g_0}^n$. We show in this
section that loops $Q(G,*,g_0,\alpha,\beta,\gamma,\delta)$ obtained with
identical $*$ are not interesting.

Let $\psi$ be a groupoid identity, and let $\var{\psi}$ be all the variables
appearing in $\psi$. Assume that for every $x\in\var{\psi}$ a decision has
been made whether $x$ is to be taken from $G$ or from $Gu$. Then, while
evaluating each side of the identity $\psi$ in $G\cup Gu$, we have to use the
multiplications $\alpha$, $\beta$, $\gamma$ and $\delta$ certain number of
times.

\begin{exam} Consider the left alternative law $x(xy)=(xx)y$. With $x\in G$,
$y\in Gu$, we see that we need $\beta$ twice to evaluate $x\circ (x\circ y)$,
while we need $\alpha$ once and $\beta$ once to evaluate $(x\circ x)\circ y$.
\end{exam}

A groupoid identity is said to be \emph{strictly balanced} if the same
variables appear on both sides of the identity the same number of times and
in the same order. For instance $(x(y(xz)))(yx) = ((xy)x)(z(yx))$ is strictly
balanced.

The above example shows that the same multiplications do not have to be used
the same number of times even while evaluating a strictly balanced identity.
However:

\begin{lemm}\lemmlbl{balanced}
Let $\psi$ be a strictly balanced identity. Assume that for $x\in\var{\psi}$
a decision has been made whether $x\in G$ or $x\in Gu$. Then, while
evaluating $\psi$ in $Q(G,*,g_0,\alpha,\beta,\gamma,\delta)$, $\delta$ is
used the same number of times on both sides of $\psi$.
\end{lemm}
\begin{proof}
Let $k$ be the number of variables on each side of $\psi$, with repetitions,
whose value is assigned to be in $Gu$. The number $k$ is well-defined since
$\psi$ is strictly balanced.

While evaluating the identity $\psi$, each multiplication reduces the number
of factors by $1$. However, only $\delta$ reduces the number of factors from
$Gu$ (by two). Since the coset multiplication in $G\cup Gu$ modulo $G$ is
associative, and since $\psi$ is strictly balanced, either both evaluated
sides of $\psi$ will end up in $G$ (in which case $\delta$ is applied $k/2$
times on each side), or both evaluated sides of $\psi$ will end up in $Gu$
(in which case $\delta$ is applied $(k-1)/2$ times on each side).
\end{proof}

\begin{lemm}\lemmlbl{neutral} If $\alpha\in\Theta$ and
$L=Q(G,*,g_0,\alpha,\beta,\gamma,\delta)$ is a loop, then the neutral element
of $Q$ coincides with the neutral element of $G$.
\end{lemm}
\begin{proof}
Let $e$ be the neutral element of $L$ and $1$ the neutral element of $G$.
Since $1=1^*$, we have $1\circ 1 = \Delta\alpha(1,1) = 1 = 1\circ e$, and the
result follows from the fact that $L$ is a quasigroup.
\end{proof}

\begin{prop} Assume that $g^*=g$ for every $g\in G$, and let $\alpha$,
$\beta$, $\gamma$, $\delta\in \Theta_0$. If
$L=Q(G,*,g_0,\alpha,\beta,\gamma,\delta)$ happens to be a loop, then every
strictly balanced identity holds in $L$. In particular, $L$ is an abelian
group.
\end{prop}
\begin{proof}
Since $*$ is identical on $G$, we have
$\Theta_0=\{\theta_{g_0}^n;\;n\in\mathbb Z\}$. By \assufref{loop}, we have
$\alpha=\theta_{xy}$. Then by \lemmfref{neutral}, $L$ has neutral element
$1$. Assume that $\beta = \theta_{g_0}^n$ for some $n$. Then $gu=1\circ gu =
(\Delta\beta(1,g))u = (g_0^ng)u$, which means that $n=0$. Similarly, if
$\gamma = \theta_{g_0}^m$ then $m=0$.

Let $\delta=\theta_{g_0}^k$. Let $\psi$ be a strictly balanced identity. For
every $x\in\var{\psi}$, decide if $x\in G$ or $x\in Gu$. By
\lemmfref{balanced}, while evaluating $\psi$ in $L$, the multiplication
$\delta$ is used the same number of times on the left and on the right, say
$t$ times. Since $\alpha=\beta=\gamma=\theta_{xy}$, we conclude that $\psi$
reduces to $g_0^{kt}z = g_0^{kt}z$, for some $z\in G\cup Gu$.

Since the associative law is strictly balanced, $L$ is associative. We have
already noticed that identical $*$ forces $G$ to be abelian. Then $L$ is
abelian too, as $gu\circ h = (gh)u = (hg)u = h\circ gu$ and $gu\circ hu =
g_0^kgh = g_0^khg = hu\circ gu$ for every $g$, $h\in G$.
\end{proof}

We have just seen that if $g=g^*$ for every $g\in G$ then our constructions
do not yield nonassociative loops. Therefore:

\begin{assu}\assulbl{nonidentical}
From now on, we assume that there exists $g\in G$ such that $g^*\ne g$.
\end{assu}

\section{Loops}\label{Sc:Loops}

In this section we further narrow the choices of $\alpha$, $\beta$, $\gamma$,
$\delta$ when $Q(G,\alpha,\beta,\gamma,\delta)$ is supposed to be a loop.

\begin{prop}\proplbl{loop}
Let $L = Q(G,*,g_0,\alpha,\beta,\gamma,\delta)$. Then $L$ is a loop if and
only if $\alpha=\theta_{xy}$, $\beta\in\{\theta_{xy}$, $\theta_{x^*y}$,
$\theta_{yx}$, $\theta_{yx^*}\}$, $\gamma\in\{\theta_{xy}$, $\theta_{xy^*}$,
$\theta_{yx}$, $\theta_{y^*x}\}$, and $\delta$ is of the form
$\theta\theta_{g_0}^n$ for some integer $n$ and $g_0\in G$.
\end{prop}
\begin{proof}
If $L$ is a loop then $\alpha\in\{\theta_{xy}$, $\theta_{xy^*}$,
$\theta_{x^*y}$, $\theta_{x^*y^*}\}$ and \lemmfref{neutral} imply that $1$ is
the neutral element of $L$.

The equation $g=1\circ g$ holds for every $g\in G$ if and only if
$\Delta\alpha(1,g)=g$ for every $g\in G$, which happens if and only if
$\alpha\in\{\theta_{xy}$, $\theta_{x^*y}\}$. (Note that we use
\assufref{nonidentical} here.) Similarly, $g=g\circ 1$ holds for every $g\in
G$ if and only if $\Delta\alpha(g,1)=g$ for every $g\in G$, which happens if
and only if $\alpha\in\{\theta_{xy}$, $\theta_{xy^*}\}$. Therefore $g=1\circ
g=g\circ 1$ holds for every $g\in G$ if and only if $\alpha=\theta_{xy}$.

Now, $gu=1\circ gu$ holds for every $g\in G$ if and only if
$\Delta\beta(1,g)=g$ for every $g\in G$, which happens if and only if
$\beta\in\{\theta_{xy}$, $\theta_{x^*y}$, $\theta_{yx}$, $\theta_{yx^*}\}$.
Similarly, $gu=gu\circ 1$ holds for every $g\in G$ if and only if
$\Delta\gamma(g,1)=g$ for every $g\in G$, which happens if and only if
$\gamma\in\{\theta_{xy}$, $\theta_{xy^*}$, $\theta_{yx}$, $\theta_{y^*x}\}$.
\end{proof}

We are only interested in loops, and we have already noted that
$(g_0^n)^*=g_0^n\in Z(G)$. Since we allow $g_0=1$, we can agree on:

\begin{assu}\assulbl{128} From now on, we assume that $\alpha =
\theta_{xy}$, $\beta\in\{\theta_{xy}$, $\theta_{x^*y}$, $\theta_{yx}$,
$\theta_{yx^*}\}$, $\gamma\in\{\theta_{xy}$, $\theta_{xy^*}$, $\theta_{yx}$,
$\theta_{y^*x}\}$, and $\delta\in\theta_{g_0}\Theta$.
\end{assu}

Our last reduction concerns the maps $\beta$ and $\gamma$.

\begin{lemm}\lemmlbl{commstar} We have
$\Delta\theta_{x^*y^*}\theta_0=\Delta\theta_0\theta_{x^*y^*}$ for every
$\theta_0\in\Theta_0$.
\end{lemm}
\begin{proof} The group $\Theta_0$ is generated by $\theta_{yx}$,
$\theta_{xy^*}$ and $\theta_{g_0}$. It therefore suffices to check that
$\Delta\theta_{x^*y^*}\theta_0=\Delta\theta_0\theta_{x^*y^*}$ holds for
$\theta_0\in\{\theta_{yx}$, $\theta_{xy^*}$, $\theta_{g_0}\}$, which follows
by straightforward calculation.
\end{proof}

\begin{lemm}\lemmlbl{iso} The quasigroups
$Q(G,*,g_0,\alpha,\beta,\gamma,\delta)$,
$Q(G,*,g_0,\alpha,\beta',\gamma',\theta_{x^*y^*}\delta)$ are isomorphic if
\begin{eqnarray*}
    (\beta,\beta')&\in&\{(\theta_{xy},\theta_{yx^*}),(\theta_{yx},\theta_{x^*y}),
        (\theta_{x^*y},\theta_{yx}), (\theta_{yx^*},\theta_{xy})\},\\
    (\gamma,\gamma')&\in&\{(\theta_{xy},\theta_{y^*x}),(\theta_{yx},\theta_{xy^*}),
        (\theta_{xy^*},\theta_{yx}),(\theta_{y^*x},\theta_{xy})\}.
\end{eqnarray*}
\end{lemm}
\begin{proof}
Let $\circ$ denote the multiplication in $Q(G,*,g_0,
\alpha,\beta,\gamma,\delta)$, and $\bullet$ the multiplication in
$Q(G,*,g_0,\alpha,\beta',\gamma',\theta_{x^*y^*}\delta)$. Consider the
permutation $f$ of $G$ defined by $f(g)=g$, $f(gu)=g^*u$, for $g\in G$.

We show that $f$ is an isomorphism of $(G\cup Gu,\circ)$ onto $(G\cup
Gu,\bullet)$ if and only if
\begin{equation}\label{Eq:Aux1}
    (\Delta\beta(g,h))^*=\Delta\beta'(g,h^*),\quad
    (\Delta\gamma(g,h))^*=\Delta\gamma'(g^*,h).
\end{equation}
Once we establish this fact, the proof is finished by checking that the pairs
$(\beta,\beta')$, $(\gamma,\gamma')$ in the statement of the Lemma satisfy
$(\ref{Eq:Aux1})$.

Let $g$, $h\in G$. Then
\begin{eqnarray*}
    f(g\circ h)&=&f(\Delta\alpha(g,h)) = \Delta\alpha(g,h),\\
    f(g\circ hu)&=&f(\Delta\beta(g,h)u) = (\Delta\beta(g,h))^*u,\\
    f(gu\circ h)&=&f(\Delta\gamma(g,h)u) = (\Delta\gamma(g,h))^*u,\\
    f(gu\circ hu)&=&f(\Delta\delta(g,h))=\Delta\delta(g,h),
\end{eqnarray*}
while
\begin{eqnarray*}
    f(g)\bullet f(h)&=& g\bullet h = \Delta\alpha(g,h),\\
    f(g)\bullet f(hu)&=& g\bullet h^*u = \Delta\beta'(g,h^*)u,\\
    f(gu)\bullet f(h)&=& g^*u\bullet h = \Delta\gamma'(g^*,h)u,\\
    f(gu)\bullet f(hu)&=& g^*u\bullet h^*u =
        \Delta\theta_{g^*h^*}\delta(g^*,h^*).
\end{eqnarray*}
We see that $f(g\circ h)=f(g)\bullet f(h)$ always holds. By
\lemmfref{commstar}, $f(gu\circ hu)=f(gu)\bullet f(hu)$ always holds.
Finally, $f(g\circ hu)=f(g)\bullet f(hu)$, $f(gu\circ h)=f(gu)\bullet f(h)$
hold if and only if $(\beta,\beta')$, $(\gamma,\gamma')$ satisfy
$(\ref{Eq:Aux1})$.
\end{proof}

Assume that $Q(G,*,g_0,\theta_{xy},\beta,\gamma,\delta)$ is a loop
(satisfying \assufref{128}). Then \lemmfref{iso} provides an isomorphism of
$Q(G,*,g_0,\theta_{xy},\beta,\gamma,\delta)$ onto some loop\linebreak
$Q(G,*,g_0,\theta_{xy},\beta',\gamma',\delta')$ such that if
$\gamma=\theta_{xy^*}$ then $\gamma' = \theta_{yx}$, and if
$\gamma=\theta_{y^*x}$ then $\gamma'=\theta_{xy}$. We can therefore assume:

\begin{assu}\assulbl{final}
From now on, we assume that $\alpha = \theta_{xy}$, $\beta\in\{\theta_{xy}$,
$\theta_{x^*y}$, $\theta_{yx}$, $\theta_{yx^*}\}$, $\gamma\in\{\theta_{xy}$,
$\theta_{yx}\}$, and $\delta\in\theta_{g_0}\Theta$.
\end{assu}

It is easy to see how much calculation is needed to find all loops
$Q(G,*,g_0,\alpha,\beta,\gamma,\delta)$ that satisfy a given groupoid
identity $\psi$. We have $1\cdot 4\cdot 2\cdot 8=64$ choices for
$(\alpha,\beta,\gamma,\delta)$. (To appreciate the reductions, compare this
with the unrestricted case $\alpha$, $\beta$, $\gamma$, $\delta\in
\Theta_0$.) Once $(\alpha,\beta,\gamma,\delta)$ is chosen, we must verify
$2^k$ equations in $G$, where $k$ is the number of variables in $\psi$ (since
each variable can be assigned value in $G$ or in $Gu$).

We work out the calculation for one identity $\psi$ and one choice of
multiplication $(\alpha,\beta,\gamma,\delta)$. After seing the routine nature
of the calculations, we gladly switch to a computer search.

\section{C-loops arising from the construction of de Barros and Juriaans}
\label{Sc:dBJ}

\noindent C-loops are loops satisfying the identity $((xy)y)z = x(y(yz))$. In
\cite{dBJ}, de Barros and Juriaans used a construction similar to
\eqref{Eq:Chein} to obtain loops whose loop algebras are flexible. In our
systematic notation, their construction is
\begin{equation}\label{Eq:dBJ}
    Q(G,*,g_0,\theta_{xy},\theta_{xy},\theta_{y^*x},\theta_{g_0}\theta_{xy^*}),
\end{equation}
with the usual conventions on $g_0$ and $*$. The construction \eqref{Eq:dBJ}
violates \assufref{final} but, by \lemmfref{iso}, it is isomorphic to
\begin{displaymath}
    Q(G,*,g_0,\theta_{xy},\theta_{yx^*},\theta_{xy},\theta_{g_0}\theta_{x^*y}),
\end{displaymath}
which complies with all assumptions we have made.

\begin{theo} Let $G$ be a group and let $L$ be the loop defined by
$(\ref{Eq:dBJ})$. Then $L$ is a flexible loop, and the following conditions
are equivalent:
\begin{enumerate}
\item[(i)] $L$ is associative,

\item[(ii)] $L$ is Moufang,

\item[(iii)] $G$ is commutative.
\end{enumerate}
Furthermore, $L$ is a C-loop if and only if $G/Z(G)$ is an elementary abelian
$2$-group. When $L$ is a C-loop, it is diassociative.
\end{theo}
\begin{proof} Throughout the proof, we use $g_0=g_0^*\in Z(G)$, $gg^*=g^*g\in
Z(G)$, $(g^*)^*=g$ and $(gh)^*=h^*g^*$ without warning.

By \propfref{loop}, $L$ is a loop.

Flexibility. For $x$, $y\in G$ we have:
\begin{align*}
    &(x\circ y)\circ x = (xy)x = x(yx) = x\circ(y\circ x),\\
    &(x\circ yu)\circ x = (xy)u\circ x = x^*xyu = xx^*yu = x\circ x^*yu =
    x\circ(yu\circ x),\\
    &(xu\circ y)\circ xu = y^*xu\circ xu = g_0y^*xx^* = g_0xx^*y^* =
    xu\circ (yx)u  = xu\circ (y\circ xu),\\
    &(xu\circ yu)\circ xu = g_0xy^*\circ xu = g_0xy^*xu =
    xu\circ g_0yx^* = xu\circ (yu\circ xu).
\end{align*}
Thus $L$ is flexible.

Associativity. For $x$, $y$, $z\in G$ we have:
\begin{align*}
    &x\circ(y\circ z) = x(yz) = (xy)z = (x\circ y)\circ z,\\
    &x\circ(y\circ zu) = x(yz)u = (xy)zu = (x\circ y)\circ zu,\\
    &xu\circ(y\circ z) = xu\circ yz = z^*y^*xu = y^*xu\circ z
    = (xu\circ y)\circ z,\\
    &x\circ(yu\circ zu) = x\circ g_0yz^* = g_0xyz^* =
    xyu\circ zu = (x\circ yu)\circ zu,\\
    &xu\circ(yu\circ z) = xu\circ z^*yu = g_0xy^*z =
    g_0xy^*\circ z= (xu\circ yu)\circ z.
\end{align*}
Furthermore,
\begin{align*}
    &x\circ(yu\circ z ) = x\circ z^*yu = xz^*yu,
    \quad\quad(x\circ yu)\circ z = xyu\circ z = z^*xyu,\\
    &xu\circ(y\circ zu) = xu\circ yzu = g_0xz^*y^*,
    \quad\quad(xu\circ y)\circ zu = y^*xu\circ zu = g_0y^*xz^*,\\
    &xu\circ(yu\circ zu) = xu\circ g_0yz^* = g_0zy^*xu,
    \quad\quad(xu\circ yu)\circ zu = g_0xy^*\circ zu = g_0xy^*zu.
\end{align*}
Thus $L$ is associative if and only if $G$ is commutative. (Sufficiency is
obvious. For necessity, note that $*$ is onto, and substitute $1$ for one of
$x$, $y$, $z$ if needed.)

Moufang property. Let $x$, $y$, $z\in G$. Then
\begin{align*}
    x\circ(yu\circ (x\circ z)) = x\circ(yu\circ xz) = x\circ z^*x^*yu = xz^*x^*yu,\\
    ((x\circ yu)\circ x)\circ z = (xyu\circ x)\circ z = x^*xyu\circ z =
    z^*x^*xyu.
\end{align*}
Therefore, this particular form of the Moufang identity holds if and only if
$xz^*x^* = z^*x^*x$. Now, given $x$, $y\in G$, there is $z\in G$ such that
$z^*x^*=y$. Therefore $xz^*x^*=z^*x^*x$ holds in $G$ if and only if $G$ is
commutative. However, when $G$ is commutative, then $L$ is associative, and
we have proved the equivalence of (i), (ii), (iii).

C property. Let $x$, $y$, $z\in G$. Then
\begin{align*}
    &x\circ(y\circ(y\circ z)) = x(y(yz)) = ((xy)y)z =
    ((x\circ y)\circ y)\circ z,\\
    &x\circ(y\circ(y\circ zu)) = (x(y(yz))u = ((xy)y)z)u =
    ((x\circ y)\circ y)\circ zu,\\
    &x\circ(yu\circ(yu\circ z)) = x\circ(yu\circ z^*yu) =
    x\circ g_0yy^*z = g_0xyy^*z=g_0xyy^*\circ z\\
    &\quad\quad= (xyu\circ yu)\circ z= ((x\circ yu)\circ yu)\circ z,\\
    &xu\circ(y\circ(y\circ z)) = xu\circ yyz = z^*y^*y^*xu
    = y^*y^*xu\circ z = (y^*xu\circ y)\circ z\\
    &\quad\quad= ((xu\circ y)\circ y)\circ z,\\
    &x\circ(yu\circ(yu\circ zu)) = x\circ (yu\circ g_0yz^*) =
    x\circ g_0zy^*yu = g_0xzy^*yu \\
    &\quad\quad = g_0xyy^*zu  = g_0xyy^*\circ zu =
    (xyu\circ yu)\circ zu  = ((x\circ yu)\circ yu)\circ zu,\\
    & xu\circ(yu\circ(yu\circ z)) = xu\circ(yu\circ z^*yu) =
    xu\circ g_0yy^*z = g_0z^*yy^*xu = g_0z^*xy^*yu\\
    &\quad\quad  = g_0xy^*yu\circ z = (g_0xy^*\circ yu)\circ z
    = ((xu\circ yu)\circ yu)\circ z,\\
    &xu\circ(yu\circ(yu\circ zu)) = xu\circ( yu\circ g_0yz^* )
    = xu\circ g_0zy^*yu = g_0^2xy^*yz^* = g_0xy^*yu\circ zu\\
    &\quad\quad = (g_0xy^*\circ yu)\circ zu = ((xu\circ yu)\circ yu)\circ zu.
\end{align*}
While verifying the remaining form of the C identity, we obtain
\begin{align*}
    &xu\circ(y\circ(y\circ zu)) = xu\circ yyzu = g_0xz^*y^*y^*,\\
    &((xu\circ y)\circ y)\circ zu = (y^*xu\circ y)\circ zu
    = y^*y^*xu\circ zu = g_0y^*y^*xz^*.
\end{align*}
The identity therefore holds if and only if $y^*y^*$ commutes with all
elements of $G$, which happens if and only if $G/Z(G)$ is an elementary
abelian $2$-group.

Finally, by Lemma 4.4 of \cite{CLoopsI}, flexible C-loops are diassociative.
\end{proof}

\section{The Algorithm}\label{Sc:Search}

\subsection{Collecting Identities}

Let $G$ be a group, $\psi$ a groupoid identity and
$(\alpha,\beta,\gamma,\delta)$ a multiplication. Then the following algorithm
will output a set $\Psi$ of group identities such that
$Q(G,*,g_0,\alpha,\beta,\gamma,\delta)$ satisfies $\psi$ if and only if $G$
satisfies all identities of $\Psi$:
\begin{enumerate}
\item[(i)] Let $f:\var{\psi}\to\{0,1\}$ be a function that decides whether
$x\in\var{\psi}$ is to be taken from $G$ or from $Gu$.

\item[(ii)] Upon assigning the variables of $\psi$ according to $f$, let
$\psi_f=(u,v)$ be the identity $\psi$ evaluated in
$Q(G,*,g_0,\alpha,\beta,\gamma,\delta)$.

\item[(iii)] Let $\Psi = \{\psi_f;\;f:\var{\psi}\to\{0,1\}\}$.
\end{enumerate}
This algorithm is straightforward but not very useful, since it typically
outputs a large number of complicated group identities.

\subsection{Understanding the identities in the Bol-Moufang case}

We managed to decipher the meaning of $\Psi$ for all multiplications
$(\alpha,\beta,\gamma,\delta)$ and for all identities of Bol-Moufang type by
another algotihm. First, we reduced the identity $\psi_f=(u,v)$ to a
canonical form as follows:
\begin{enumerate}
\item[(a)] replace $g_0^*$ by $g_0$,

\item[(b)] move all $g_0$ to the very left,

\item[(c)] replace $x^*x$ by $xx^*$,

\item[(d)] move all substrings $xx^*$ immediately to the right of the power
$g_0^m$, and order the substrings $xx^*$, $yy^*$, $\dots$ lexicographically,

\item[(e)] cancel as much as possible on the left and on the right of the
resulting identity.
\end{enumerate}
Then we used Lemmas \lemmref{lm1}--\lemmref{lm5} to understand what the
canonical identities collected in $\Psi$ say about the group $G$:

\begin{lemm}\lemmlbl{lm1}
If an identity of $\Psi$ reduces to $x^*=x$ then it does not hold in any
group.
\end{lemm}
\begin{proof}
Since we assume that $*$ is not identical on $G$.
\end{proof}

\begin{lemm} The following conditions are equivalent:
\begin{enumerate}
\item[(i)] $G/Z(G)$ is an elementary abelian $2$-group,

\item[(ii)] $xxy=yxx$,

\item[(iii)] $xyx^*=x^*yx$.
\end{enumerate}
\end{lemm}
\begin{proof}
We have $xyx^*=x^*yx$ if and only if $x^*xyx^*x = x^*x^*yxx$. Since $x^*x\in
Z(G)$, the latter identity is equivalent to $x^*xx^*xy = x^*x^*yxx$. Since
$xx^*=x^*x$, we can rewrite it equivalently as $x^*x^*xxy=x^*x^*yxx$, which
is by cancellation equivalent to $xxy=yxx$.
\end{proof}

\begin{lemm} The following conditions are equivalent:
\begin{enumerate}
\item[(i)] $G$ is commutative,

\item[(ii)] $xx^*y=x^*yx$.
\end{enumerate}
\end{lemm}
\begin{proof}
If $xx^*y=x^*yx$ then $x^*xy=x^*yx$ and so $xy=yx$.
\end{proof}

\begin{lemm} If $\psi$ is a strictly balanced identity that reduces to
$xy=yx$ upon substituting $1$ for some of the variables of $\psi$, then
$\psi$ is equivalent to commutativity.
\end{lemm}
\begin{proof} $\psi$ implies commutativity. Once commutativity holds, we
can rearrange the variables of $\psi$ so that both sides of $\psi$ are the
same, because $\psi$ is strictly balanced.
\end{proof}

\begin{lemm}\lemmlbl{lm5} The following conditions are equivalent:
\begin{enumerate}
\item[(i)] $xxy=yx^*x^*$ holds in $G$,

\item[(ii)] $(xx)^*=xx$ and $G/Z(G)$ is an elementary abelian $2$-group.
\end{enumerate}
\end{lemm}
\begin{proof}
Condition (ii) clearly implies (i). If (i) holds, we have $xx=x^*x^*$ (with
$y=1$) and so $(xx)^*=xx$. Also $xxy=yx^*x^*=yxx$.
\end{proof}

\subsection{What the identities mean in the Bol-Moufang case}

Lemmas \lemmref{lm1}--\lemmref{lm5} are carefully tailored to loops of
Bol-Moufang type, and we discovered them upon studying the canonical
identities $\Psi$ obtained by the computer search.

It just so happens that every identity $\psi_f$ of $\Psi$ is equivalent to a
combination of the following properties of $G$:
\begin{enumerate}
\item[(PN)] No group satisfies $\psi_f$.

\item[(PA)] All groups satisfy $\psi_f$.

\item[(PC)] $G$ is commutative.

\item[(PB)] $G/Z(G)$ is an elementary abelian $2$-group.

\item[(PS)] $(gg)^*=gg$ for every $g\in G$.
\end{enumerate}
A prominent example of $*$ is the inverse operation ${}^{-1}$ in $G$. Then
(PB) says that $G$ is of exponent $4$, and it is therefore not difficult to
obtain examples of groups satisfying any possible combination of (PN), (PA),
(PC), (PB) and (PS).

We have implemented the algorithm in GAP \cite{GAP}, and made it available
online at
\begin{displaymath}
\texttt{http://www.math.du.edu/\symbol{126}petr}
\end{displaymath}
in section Research. The algorithm is not safe for identities that are not
strictly balanced.

\section{Results}\label{Sc:Results}

We now present the results of the computer search. In order to organize the
results, observe that if $L=Q(G,*,g_0,\alpha,\beta,\gamma,\delta)$ is
associative, it satisfies all identities of Bol-Moufang type. Since we do not
want to list the multiplications and properties of $G$ repeatedly, we first
describe all cases when $L$ is associative, then all cases when $L$ is an
extra loop, then all cases when $L$ is a Moufang loop, etc., guided by the
inclusions of Figure \ref{Fg:BM}.

All results of this section are computer generated. To avoid errors in
transcribing, the \TeX\ source of the statements of the results is also
computer generated. In the statements, we write $xy$ instead of
$\theta_{xy}$, $g_0yx^*$ instead of $\theta_{g_0}\theta_{yx^*}$, etc., in
order to save space and improve legibility. Some results are mirror versions
of others (cf. \theofref{leftBol} versus \theofref{rightBol}), but we decided
to include them anyway for quicker future reference. Finally, when $G$ is
commutative, $\Delta(\Theta\cup\theta_{g_0}\Theta)$ coincides with
$\Delta(S\cup \theta_{g_0}S)$, where $S=\{\theta_{xy}$, $\theta_{xy^*}$,
$\theta_{x^*y}$, $\theta_{x^*y^*}\}$. We therefore report only maps $\alpha$,
$\beta$, $\gamma$, $\delta$ from $S\cup \theta_{g_0}S$ in the commutative
case.

In Theorems \theoref{associative} -- \theoref{rns}, $G$ is a group, $*$ is a
nonidentical involutory antiautomorphism of $G$ satisfying $gg^*\in Z(G)$ for
every $g\in G$, the element $g_0\in Z(G)$ satisfies $g_0^*=g_0$, and the maps
$\alpha$, $\beta$, $\gamma$, $\delta$ are as in \assufref{final}.

\begin{theo}\theolbl{associative}
The loop $Q(G,*,g_0,\theta_{xy},\beta,\gamma,\delta)$ is
associative iff the following conditions are satisfied: \vskip
2mm\noindent$(\beta,\gamma,\delta)$ is
equal to\\
$(xy,xy,g_0xy)$, or

\vskip 2mm\noindent$G$ is commutative and $(\beta,\gamma,\delta)$
is equal to $(x^*y,xy,g_0x^*y)$.
% manually removed $(xy,xy,g_0xy)$,
\end{theo}

\begin{theo}
The loop $Q(G,*,g_0,\theta_{xy},\beta,\gamma,\delta)$ is extra iff
it is associative or if the following conditions are satisfied:

\vskip 2mm\noindent$G/Z(G)$ is an
 elementary abelian $2$-group
 and $(\beta,\gamma,\delta)$ is
equal to\\
$(x^*y,yx,g_0yx^*)$.
\end{theo}

\begin{theo}\theolbl{Moufang}
The loop $Q(G,*,g_0,\theta_{xy},\beta,\gamma,\delta)$ is Moufang
iff it is extra or if the following conditions are satisfied:
\vskip 2mm\noindent$(\beta,\gamma,\delta)$ is
equal to\\
$(x^*y,yx,g_0yx^*)$.
\end{theo}

\begin{theo}
The loop $Q(G,*,g_0,\theta_{xy},\beta,\gamma,\delta)$ is a C-loop
iff it is extra or if the following conditions are satisfied:

\vskip 2mm\noindent$G/Z(G)$ is an
 elementary abelian $2$-group
 and $(\beta,\gamma,\delta)$ is
among\\
$(yx,yx,g_0yx)$, $(yx^*,xy,g_0x^*y)$.
\end{theo}

\begin{theo}\theolbl{leftBol}
The loop $Q(G,*,g_0,\theta_{xy},\beta,\gamma,\delta)$ is left Bol
iff it is Moufang or if the following conditions are satisfied:

\vskip 2mm\noindent$G/Z(G)$ is an
 elementary abelian $2$-group
 and $(\beta,\gamma,\delta)$ is
among\\
$(xy,yx,g_0yx)$, $(x^*y,xy,g_0x^*y)$, or

\vskip 2mm\noindent$G$ is commutative, $(xx)^*=xx$ for
 every $x\in G$ and $(\beta,\gamma,\delta)$ is
among\\
$(xy,xy,g_0x^*y)$, $(x^*y,xy,g_0xy)$, or

\vskip 2mm\noindent$G/Z(G)$ is
 an elementary abelian $2$-group,
 $(xx)^*=xx$ for
 every $x\in G$ and $(\beta,\gamma,\delta)$ is
among\\
$(xy,xy,g_0x^*y)$, $(xy,yx,g_0yx^*)$, $(x^*y,xy,g_0xy)$,
$(x^*y,yx,g_0yx)$.
\end{theo}

\begin{theo}\theolbl{rightBol}
The loop $Q(G,*,g_0,\theta_{xy},\beta,\gamma,\delta)$ is right Bol
iff it is Moufang or if the following conditions are satisfied:

\vskip 2mm\noindent$G/Z(G)$ is an
 elementary abelian $2$-group
 and $(\beta,\gamma,\delta)$ is
among\\
$(yx,xy,g_0yx)$, $(yx^*,yx,g_0x^*y)$, or

\vskip 2mm\noindent$G$ is commutative, $(xx)^*=xx$ for
 every $x\in G$ and $(\beta,\gamma,\delta)$ is
among\\
$(xy,xy,g_0xy^*)$, $(x^*y,xy,g_0x^*y^*)$, or

\vskip 2mm\noindent$G/Z(G)$ is
 an elementary abelian $2$-group,
 $(xx)^*=xx$ for
 every $x\in G$ and $(\beta,\gamma,\delta)$ is
among\\
$(xy,xy,g_0xy^*)$, $(x^*y,yx,g_0y^*x^*)$, $(yx,xy,g_0y^*x)$,
$(yx^*,yx,g_0x^*y^*)$.
\end{theo}

\begin{theo}
The loop $Q(G,*,g_0,\theta_{xy},\beta,\gamma,\delta)$ is an
LC-loop iff it is a C-loop or if the following conditions are
satisfied:

\vskip 2mm\noindent$G/Z(G)$ is an
 elementary abelian $2$-group
 and $(\beta,\gamma,\delta)$ is
among\\
$(xy,yx,g_0yx)$, $(x^*y,xy,g_0x^*y)$, $(yx,xy,g_0xy)$,
$(yx^*,yx,g_0yx^*)$, or

\vskip 2mm\noindent$G$ is commutative, $(xx)^*=xx$ for
 every $x\in G$ and $(\beta,\gamma,\delta)$ is
among\\
$(xy,xy,g_0x^*y)$, $(x^*y,xy,g_0xy)$, or

\vskip 2mm\noindent$G/Z(G)$ is
 an elementary abelian $2$-group,
 $(xx)^*=xx$ for
 every $x\in G$ and $(\beta,\gamma,\delta)$ is
among\\
$(xy,xy,g_0x^*y)$, $(xy,yx,g_0yx^*)$, $(x^*y,xy,g_0xy)$,
$(x^*y,yx,g_0yx)$,\\
$(yx,xy,g_0x^*y)$, $(yx,yx,g_0yx^*)$, $(yx^*,xy,g_0xy)$,
$(yx^*,yx,g_0yx)$.
\end{theo}

\begin{theo}
The loop $Q(G,*,g_0,\theta_{xy},\beta,\gamma,\delta)$ is an
RC-loop iff it is a C-loop or if the following conditions are
satisfied:

\vskip 2mm\noindent$G/Z(G)$ is an
 elementary abelian $2$-group
 and $(\beta,\gamma,\delta)$ is
among\\
$(xy,yx,g_0xy)$, $(x^*y,xy,g_0yx^*)$, $(yx,xy,g_0yx)$,
$(yx^*,yx,g_0x^*y)$, or

\vskip 2mm\noindent$G$ is commutative, $(xx)^*=xx$ for
 every $x\in G$ and $(\beta,\gamma,\delta)$ is
among\\
$(xy,xy,g_0xy^*)$, $(x^*y,xy,g_0x^*y^*)$, or

\vskip 2mm\noindent$G/Z(G)$ is
 an elementary abelian $2$-group,
 $(xx)^*=xx$ for
 every $x\in G$ and $(\beta,\gamma,\delta)$ is
among\\
$(xy,xy,g_0xy^*)$, $(xy,yx,g_0xy^*)$, $(x^*y,xy,g_0y^*x^*)$,
$(x^*y,yx,g_0y^*x^*)$,\\
$(yx,xy,g_0y^*x)$, $(yx,yx,g_0y^*x)$, $(yx^*,xy,g_0x^*y^*)$,
$(yx^*,yx,g_0x^*y^*)$.
\end{theo}

\begin{theo}
The loop $Q(G,*,g_0,\theta_{xy},\beta,\gamma,\delta)$ is flexible
iff it is Moufang or if the following conditions are satisfied:
\vskip 2mm\noindent$(\beta,\gamma,\delta)$ is
among\\
$(xy,xy,g_0y^*x^*)$, $(x^*y,yx,g_0xy^*)$, $(x^*y,yx,g_0x^*y)$,
$(x^*y,yx,g_0y^*x)$,\\
$(yx,yx,g_0x^*y^*)$, $(yx,yx,g_0yx)$, $(yx^*,xy,g_0xy^*)$,
$(yx^*,xy,g_0x^*y)$,\\
$(yx^*,xy,g_0yx^*)$, $(yx^*,xy,g_0y^*x)$, or

%manually removed
%\vskip 2mm\noindent$G$ is commutative and $(\beta,\gamma,\delta)$ is
%among\\
%$(xy,xy,g_0x^*y^*)$,
%$(x^*y,xy,g_0xy^*)$, or

\vskip 2mm\noindent$G/Z(G)$ is an
 elementary abelian $2$-group
 and $(\beta,\gamma,\delta)$ is
among\\
$(xy,xy,g_0x^*y^*)$, $(xy,xy,g_0yx)$, $(yx,yx,g_0xy)$,
$(yx,yx,g_0y^*x^*)$.
\end{theo}

\begin{theo}
The loop $Q(G,*,g_0,\theta_{xy},\beta,\gamma,\delta)$ is left
alternative iff it is
 left Bol or an LC-loop or if the following conditions are satisfied:
\vskip 2mm\noindent$(\beta,\gamma,\delta)$ is
among\\
$(xy,xy,g_0x^*y)$, $(xy,yx,g_0yx^*)$, $(x^*y,xy,g_0x^*y)$,
$(yx,xy,g_0x^*y)$,\\
$(yx,yx,g_0yx^*)$, $(yx^*,xy,g_0x^*y)$, $(yx^*,yx,g_0yx^*)$, or

%manually removed
%\vskip 2mm\noindent$G$ is commutative and $(\beta,\gamma,\delta)$ is
%equal to\\
%$(xy,xy,g_0x^*y)$, or

\vskip 2mm\noindent$(xx)^*=xx$ for every $x\in G$
 and $(\beta,\gamma,\delta)$ is
equal to\\
$(x^*y,xy,g_0xy)$.
\end{theo}

\begin{theo}
The loop $Q(G,*,g_0,\theta_{xy},\beta,\gamma,\delta)$ is right
alternative iff it is
 right Bol or an RC-loop or if the following conditions are satisfied:
\vskip 2mm\noindent$(\beta,\gamma,\delta)$ is
among\\
$(xy,xy,g_0xy^*)$, $(xy,yx,g_0xy^*)$, $(x^*y,xy,g_0yx^*)$,
$(yx,xy,g_0y^*x)$,\\
$(yx,yx,g_0y^*x)$, $(yx^*,xy,g_0x^*y)$, $(yx^*,yx,g_0x^*y)$, or

%manually removed
%\vskip 2mm\noindent$G$ is commutative and $(\beta,\gamma,\delta)$ is
%equal to\\
%$(xy,xy,g_0xy^*)$, or

\vskip 2mm\noindent$(xx)^*=xx$ for every $x\in G$
 and $(\beta,\gamma,\delta)$ is
equal to\\
$(yx^*,yx,g_0x^*y^*)$.
\end{theo}

\begin{theo}
The loop $Q(G,*,g_0,\theta_{xy},\beta,\gamma,\delta)$ is a left
nuclear square loop
 iff it is an LC-loop or if the following conditions are satisfied:
\vskip 2mm\noindent$(\beta,\gamma,\delta)$ is
among\\
$(xy,xy,g_0xy^*)$, $(yx^*,yx,g_0x^*y)$, $(yx^*,yx,g_0x^*y^*)$, or

\vskip 2mm\noindent$G/Z(G)$ is an
 elementary abelian $2$-group
 and $(\beta,\gamma,\delta)$ is
among\\
$(xy,xy,g_0yx)$, $(xy,xy,g_0y^*x)$, $(xy,yx,g_0xy)$,
$(xy,yx,g_0xy^*)$,\\
$(xy,yx,g_0y^*x)$, $(x^*y,xy,g_0x^*y^*)$, $(x^*y,xy,g_0yx^*)$,
$(x^*y,xy,g_0y^*x^*)$,\\
$(x^*y,yx,g_0x^*y)$, $(x^*y,yx,g_0x^*y^*)$, $(x^*y,yx,g_0y^*x^*)$,
$(yx,xy,g_0xy^*)$,\\
$(yx,xy,g_0yx)$, $(yx,xy,g_0y^*x)$, $(yx,yx,g_0xy)$,
$(yx,yx,g_0xy^*)$,\\
$(yx,yx,g_0y^*x)$, $(yx^*,xy,g_0x^*y^*)$, $(yx^*,xy,g_0yx^*)$,
$(yx^*,xy,g_0y^*x^*)$,\\
$(yx^*,yx,g_0y^*x^*)$, or

\vskip 2mm\noindent$G/Z(G)$ is
 an elementary abelian $2$-group,
 $(xx)^*=xx$ for
 every $x\in G$ and $(\beta,\gamma,\delta)$ is
among\\
$(xy,xy,g_0x^*y^*)$, $(xy,xy,g_0yx^*)$, $(xy,xy,g_0y^*x^*)$,
$(xy,yx,g_0x^*y)$,\\
$(xy,yx,g_0x^*y^*)$, $(xy,yx,g_0y^*x^*)$, $(x^*y,xy,g_0xy^*)$,
$(x^*y,xy,g_0yx)$,\\
$(x^*y,xy,g_0y^*x)$, $(x^*y,yx,g_0xy)$, $(x^*y,yx,g_0xy^*)$,
$(x^*y,yx,g_0y^*x)$,\\
$(yx,xy,g_0x^*y^*)$, $(yx,xy,g_0yx^*)$, $(yx,xy,g_0y^*x^*)$,
$(yx,yx,g_0x^*y)$,\\
$(yx,yx,g_0x^*y^*)$, $(yx,yx,g_0y^*x^*)$, $(yx^*,xy,g_0xy^*)$,
$(yx^*,xy,g_0yx)$,\\
$(yx^*,xy,g_0y^*x)$, $(yx^*,yx,g_0xy)$, $(yx^*,yx,g_0xy^*)$,
$(yx^*,yx,g_0y^*x)$.
\end{theo}

\begin{theo}
The loop $Q(G,*,g_0,\theta_{xy},\beta,\gamma,\delta)$ is a middle
nuclear square loop
 iff it is an LC-loop or an RC-loop or if the following conditions are
 satisfied:
\vskip 2mm\noindent$(\beta,\gamma,\delta)$ is
among\\
$(xy,xy,g_0y^*x^*)$, $(yx^*,xy,g_0xy^*)$, $(yx^*,xy,g_0yx^*)$, or

\vskip 2mm\noindent$G/Z(G)$ is an
 elementary abelian $2$-group
 and $(\beta,\gamma,\delta)$ is
among\\
$(xy,xy,g_0x^*y^*)$, $(xy,xy,g_0yx)$, $(xy,yx,g_0x^*y^*)$,
$(xy,yx,g_0y^*x^*)$,\\
$(x^*y,xy,g_0xy^*)$, $(x^*y,xy,g_0y^*x)$, $(x^*y,yx,g_0xy^*)$,
$(x^*y,yx,g_0x^*y)$,\\
$(x^*y,yx,g_0y^*x)$, $(yx,xy,g_0x^*y^*)$, $(yx,xy,g_0y^*x^*)$,
$(yx,yx,g_0xy)$,\\
$(yx,yx,g_0x^*y^*)$, $(yx,yx,g_0y^*x^*)$, $(yx^*,xy,g_0y^*x)$,
$(yx^*,yx,g_0xy^*)$,\\
$(yx^*,yx,g_0y^*x)$, or

\vskip 2mm\noindent$G/Z(G)$ is
 an elementary abelian $2$-group,
 $(xx)^*=xx$ for
 every $x\in G$ and $(\beta,\gamma,\delta)$ is
among\\
$(xy,xy,g_0yx^*)$, $(xy,xy,g_0y^*x)$, $(xy,yx,g_0x^*y)$,
$(xy,yx,g_0y^*x)$,\\
$(x^*y,xy,g_0x^*y^*)$, $(x^*y,xy,g_0yx)$, $(x^*y,yx,g_0xy)$,
$(x^*y,yx,g_0x^*y^*)$,\\
$(yx,xy,g_0xy^*)$, $(yx,xy,g_0yx^*)$, $(yx,yx,g_0xy^*)$,
$(yx,yx,g_0x^*y)$,\\
$(yx^*,xy,g_0yx)$, $(yx^*,xy,g_0y^*x^*)$, $(yx^*,yx,g_0xy)$,
$(yx^*,yx,g_0y^*x^*)$.
\end{theo}

\begin{theo}\theolbl{rns}
The loop $Q(G,*,g_0,\theta_{xy},\beta,\gamma,\delta)$ is a right
nuclear square loop
 iff it is an RC-loop or if the following conditions are satisfied:
\vskip 2mm\noindent$(\beta,\gamma,\delta)$ is
among\\
$(xy,xy,g_0x^*y)$, $(x^*y,xy,g_0xy)$, $(x^*y,xy,g_0x^*y)$, or

\vskip 2mm\noindent$G/Z(G)$ is an
 elementary abelian $2$-group
 and $(\beta,\gamma,\delta)$ is
among\\
$(xy,xy,g_0yx)$, $(xy,xy,g_0yx^*)$, $(xy,yx,g_0x^*y)$,
$(xy,yx,g_0yx)$,\\
$(xy,yx,g_0yx^*)$, $(x^*y,xy,g_0yx)$, $(x^*y,yx,g_0xy)$,
$(x^*y,yx,g_0x^*y)$,\\
$(x^*y,yx,g_0yx)$, $(yx,xy,g_0xy)$, $(yx,xy,g_0x^*y)$,
$(yx,xy,g_0yx^*)$,\\
$(yx,yx,g_0xy)$, $(yx,yx,g_0x^*y)$, $(yx,yx,g_0yx^*)$,
$(yx^*,xy,g_0xy)$,\\
$(yx^*,xy,g_0yx)$, $(yx^*,xy,g_0yx^*)$, $(yx^*,yx,g_0xy)$,
$(yx^*,yx,g_0yx)$,\\
$(yx^*,yx,g_0yx^*)$, or

\vskip 2mm\noindent$G/Z(G)$ is
 an elementary abelian $2$-group,
 $(xx)^*=xx$ for
 every $x\in G$ and $(\beta,\gamma,\delta)$ is
among\\
$(xy,xy,g_0x^*y^*)$, $(xy,xy,g_0y^*x)$, $(xy,xy,g_0y^*x^*)$,
$(xy,yx,g_0x^*y^*)$,\\
$(xy,yx,g_0y^*x)$, $(xy,yx,g_0y^*x^*)$, $(x^*y,xy,g_0xy^*)$,
$(x^*y,xy,g_0x^*y^*)$,\\
$(x^*y,xy,g_0y^*x)$, $(x^*y,yx,g_0xy^*)$, $(x^*y,yx,g_0x^*y^*)$,
$(x^*y,yx,g_0y^*x)$,\\
$(yx,xy,g_0xy^*)$, $(yx,xy,g_0x^*y^*)$, $(yx,xy,g_0y^*x^*)$,
$(yx,yx,g_0xy^*)$,\\
$(yx,yx,g_0x^*y^*)$, $(yx,yx,g_0y^*x^*)$, $(yx^*,xy,g_0xy^*)$,
$(yx^*,xy,g_0y^*x)$,\\
$(yx^*,xy,g_0y^*x^*)$, $(yx^*,yx,g_0xy^*)$, $(yx^*,yx,g_0y^*x)$,
$(yx^*,yx,g_0y^*x^*)$.
\end{theo}

\section{Concluding remarks}

(I) Figure \ref{Fg:BM} and Theorems \theoref{associative}--\theoref{rns}
taken together tell us more than if we consider them separately. For
instance, Theorem \theoref{associative} and Theorem \theoref{Moufang} plus
the fact that every group is a Moufang loop imply that the construction of
Theorem \theoref{Moufang} yields a nonassociative loop if and only if the
group $G$ is not commutative. In other words, the two theorems encompass
Theorem \theoref{Chein}, and, in addition, show that Chein's construction is
unique for Moufang loops.

(II) Note that we have also recovered (an isomorphic copy of) the
construction \eqref{Eq:dBJ} of de Barros and Juriaans. Our results on Bol
loops agree with those of \cite{VojtechovskyCMUC2}, obtained by hand.

(III) To illustrate how the algorithm works for loops that are not of
Bol-Moufang type, we show the output for nonassociative RIF loops. A loop is
an \emph{RIF loop} if it satisfies $(xy)(z(xy))=((x(yz))x)y$.

\begin{theo}
The loop $Q(G,*,g_0,\theta_{xy},\beta,\gamma,\delta)$ is RIF iff
it is associative or if the following conditions are satisfied:
\vskip 2mm\noindent$(\beta,\gamma,\delta)$ is
among\\
$(x^*y,yx,g_0yx^*)$, $(yx^*,xy,g_0x^*y)$, or

\vskip 2mm\noindent$(\beta,\gamma,\delta)$ and $G$ are as
 in the following list:\\
$(yx,yx,g_0yx)$ and $xyzxy=yxzyx$.
\end{theo}

Note that the algorithm did not manage to decipher the meaning of the group
identity $xyzxy=yxzyx$, so it simply listed it.

(IV) We conclude the paper with the following observation:

\begin{lemm} Let $L=Q(G,*,g_0,\alpha,\beta,\gamma,\delta)$ be a loop. Then
$L$ has two-sided inverses.
\end{lemm}
\begin{proof}
Let $g\in G$. Since $g^*(g^{-1})^* = (g^{-1}g)^*=1^*=1$, we have
$(g^*)^{-1}=(g^{-1})^*$, and the antiautomorphisms ${}^{-1}$ and $*$ commute.
Let us denote $(g^{-1})^*=(g^*)^{-1}$ by $g^{-*}$.

We show that for every $\alpha\in\Theta_0$ and $g\in G$, there is $h\in G$
such that $\Delta\alpha(g,h) = g\circ h = 1 = h\circ g = \Delta \alpha(h,g)$.
The proof for $gu\in Gu$ is similar.

Assume that $\alpha\in\{\theta_{xy}$, $\theta_{xy^*}$, $\theta_{x^*y}$,
$\theta_{x^*y^*}\}$. Then
\begin{align*}
&\Delta\theta_{xy}(g,g^{-1}) =
gg^{-1}=1=g^{-1}g=\Delta\theta_{xy}(g^{-1},g),\\
&\Delta\theta_{xy^*}(g,g^{-*}) = g(g^{-*})^* = 1 = g^{-*}g^* =
\Delta\theta_{xy^*}(g^{-*},g),\\
&\Delta\theta_{x^*y}(g,g^{-*}) = g^*g^{-*} = 1 = (g^{-*})^*g =
\Delta\theta_{x^*y}(g^{-*},g),\\
&\Delta\theta_{x^*y^*}(g,g^{-1}) = g^*g^{-*} = 1 = g^{-*}g^* =
\Delta\theta_{x^*y^*}(g^{-1},g)
\end{align*}
show that the two-sided inverse $h$ exists. The case
$\alpha\in\{\theta_{yx}$, $\theta_{y^*x}$, $\theta_{yx^*}$,
$\theta_{y^*x^*}\}$ is similar. The general case $\alpha\in\Theta_0$ then
follows thanks to $g_0 = g_0^*\in Z(G)$.
\end{proof}

\bibliographystyle{plain}

\begin{thebibliography}{10}

\bibitem[Bruck(1971)]{Bruck} R.~Hubert~Bruck, A Survey of Binary Systems, third
printing, corrected, \emph{Ergebnisse der Mathematik und ihrer Grenzgebiete},
\emph{Neue Folge} \textbf{20}, Springer-Verlag, 1971.

\bibitem[de Barros and Juriaans(1998)]{dBJ} Luis G.~X.~de~Barros and Stanley
O.~Juriaans, \emph{Some loops whose loop algebras are flexible II},
International Journal of Mathematics, Game Theory and Algebra \textbf{8}, no.
\textbf{1}, 73--80.

\bibitem[Chein(1978)]{Chein} Orin Chein, Moufang loops of small order, \emph{Mem.
Amer. Math. Soc.} \textbf{13}(1978), no. \textbf{197}.

\bibitem[Fenyves(1969)]{Fenyves} Ferenc Fenyves, \emph{Extra loops II, On loops with
identities of Bol-Moufang type}, Publ. Math. Debrecen \textbf{16}(1969),
187--192.

\bibitem[GAP (1999)]{GAP} The GAP Group, GAP --- Groups, Algorithms, and
Programming, Version 4.3; Aachen, St Andrews (1999). (Visit
http://www-gap.dcs.st-and.ac.uk/\~{}gap).

\bibitem[Pflugfelder(1990)]{Pflugfelder}
H.~O.~Pflugfelder, Quasigroups and Loops: Introduction, \emph{Sigma series in
pure mathematics} {\bf 7}, Heldermann Verlag Berlin, 1990.

\bibitem[Phillips and Vojt\v{e}chovsk\'y(2004)]{CLoopsI}
J.~D.~Phillips and Petr Vojt\v{e}chovsk\'y, \emph{On C-loops}, to appear in
Publ.\ Math.\ Debrecen.

\bibitem[Phillips and Vojt\v{e}chovsk\'y(2005)]{PhillipsVojtechovsky}
J.~D.~Phillips and Petr Vojt\v{e}chovsk\'y, \emph{The varieties of loops of
Bol-Moufang type}, to appear in Algebra Universalis.

\bibitem[Vojt\v{e}chovsk\'y(2003)]{VojtechovskyCMUC1}
Petr Vojt\v{e}chovsk\'y, \emph{On the uniqueness of loops $M(G,2)$}, Comment.
Math. Univ. Carolin. \textbf{44} (2003),  no. \textbf{4}, 629--635.

\bibitem[Vojt\v{e}chovsk\'y(2004)]{VojtechovskyCMUC2}
Petr Vojt\v{e}chovsk\'y, \emph{A class of Bol loops with a subgroup of index
two}, Comment. Math. Univ. Carolin.  \textbf{45} (2004),  no. \textbf{2},
371--381.


\end{thebibliography}

\end{document}